\documentclass[12pt,a4paper]{article}
\advance\topmargin by -1in \advance\textheight by 1.5in
\advance\oddsidemargin by -.8in \advance\textwidth by 1.7in

\author{\sf Fr\'{e}d\'{e}ric Ja\"{e}ck and Stephen Power}
\date{  }

\usepackage{amsthm}
\usepackage{amsmath}
\usepackage{amsfonts}
\usepackage{amssymb}

\newtheorem{thm}{\sf Theorem}[section]
\newtheorem{cor}[thm]{\sf Corollary}

\newtheorem{prop}[thm]{\sf Proposition}
\newenvironment{demo}{{\sf \noindent Proof}:}{$\blacksquare$\par}
\theoremstyle{definition}
\newtheorem{defn}[thm]{\sf Definition}

\numberwithin{equation}{section}

\newcommand{\Aa}{\mathcal A}
\newcommand{\Bb}{\mathcal B}
\newcommand{\CC}{\mathbb{C}}
\newcommand{\Cc}{\mathcal C}

\newcommand{\FF}{\mathbb{F}}
\newcommand{\Dd}{\mathcal D}
\newcommand{\Ee}{\mathcal E}
\newcommand{\Ff}{\mathcal F}
\newcommand{\Gg}{\mathcal G}
\newcommand{\Hh}{\mathcal H}
\newcommand{\Kk}{\mathcal K}

\newcommand{\Ll}{\mathcal L}
\newcommand{\Mm}{\mathcal M}
\newcommand{\Pp}{\mathcal P}
\newcommand{\Qq}{\mathcal Q}
\newcommand{\Ss}{\mathcal S}
\newcommand{\TT}{\mathbb{T}}
\newcommand{\Tt}{\mathcal T}

\newcommand{\Vv}{\mathcal V}
\newcommand{\Ww}{\mathcal W}
\newcommand{\Xx}{\mathcal X}
\renewcommand{\AA}{\mathbb{A}}
\newcommand{\ZZ}{\mathbb{Z}}
\renewcommand{\TT}{\mathbb{T}}

\title{Hyper-reflexivity of free semigroupoid algebras}
\begin{document}
\maketitle

\footnote{
This work is part of the research program of the network "Analysis and Operators"
supported by the European Community's Potential Programme under HPRN-CT-2000-00116
(Analysis and operators).\\
Mathematics Subject Classification 2000, 47L75. }

\begin{abstract}
As a generalization of the free semigroup algebras considered by
Davidson and Pitts, and others, the second author and D.W. Kribs
initiated a study of reflexive algebras associated with
directed graphs. A free semigroupoid algebra $\Ll_G$ is
generated by a family of partial isometries, and  initial
projections, which act on a generalized Fock
space spawned by the directed graph $G$.
We show that if the graph is finite then $\Ll_G$ is hyper-reflexive.
\end{abstract}
\section{\sf Introduction.}

The free semigroup algebra $\Ll_n$
is the weakly closed operator algebra generated by the
left regular Hilbert space representation of the free semigroup on $n$
generators. These algebras, which are free analogues of the
operator algebra $H^\infty$, have been  considered in detail by
Davidson, Pitts, Popescu and others (\cite{DKP},\cite{DP1}, \cite{DP2},
\cite{Ari-Pop}, \cite{Pop_disc}, \cite{Pop_jot99}).
 As a generalisation of this
class the second author and D.W. Kribs in \cite{KribsPower1},
\cite{KribsPower2} have considered
reflexive algebras associated with directed graphs, the so-called
free semigroupoid algebras.  Such an algebra
$\Ll_G$ is generated by a family of
partial isometries and  initial projections which act on a
generalized Fock space $\Hh_G$ spawned by a  directed
graph $G$.

Let $\Bb (\Hh)$ be the
algebra of all bounded operators on the separable Hilbert space
$\Hh$, and let $\Aa$ be a {\small WOT}-closed subspace of
$\Bb(\Hh)$. For a given operator $T\in \Bb(\Hh)$, consider the
following two quantities :
$$d(T,\Aa) :=  \inf_{A\in\Aa}\|T-A\| \mbox{ ~~and~~ }
\beta_\Aa (T) :=  \sup_{(Q,P)\in\Ff_\Aa } \|QTP\|, $$ where
$\Ff_\Aa$ is the set of all pairs of projections annihilating
$\Aa$, that is,
 $$\Ff_\Aa =  \{ (Q,P)  ;  QAP= 0,
\mbox{ for all } A \in \Aa\}.$$
Recall that $\Aa$ is said to be a reflexive space if
$\beta_\Aa(A) = 0$ entails $ A \in \Aa$. In this case $\beta_\Aa(T)$
can also be viewed as a distance from $T$ to $\Aa$.
It is easy to check that $\beta_\Aa (T) \leq d(T,\Aa)$
while if the two distances are
equivalent then the
space $\Aa$ is said to be hyper-reflexive.
 The smallest constant $C>0$ such that $d(T,\Aa)\leq
C\beta_\Aa(T)$ for every operator $T$ is known as the
hyper-reflexivity constant for $\Aa$.
This metric strengthening of reflexivity was examined by Arveson in \cite{arv-ten}.
See Davidson \cite{dav-book} for additional discussions.

There are few systematic tools available to prove
hyper-reflexivity and as a result the variety of known
hyper-reflexive algebras is limited. The hyper-reflexivity of
$\Ll_1$, which is the algebra
$H^\infty(\TT)$, realized as the T{\oe}plitz algebra on $H^2(\TT)$, was
established by K.R. Davidson  \cite{Davidson1}. By
different methods Davidson and Pitts proved that the
remaining free semigroup algebras $\Ll_n$ are
hyper-reflexive. (Reflexivity was shown earlier by Arias and Popescu
\cite{Ari-Pop} and an alternative short proof can be found in
Kribs and Power \cite{KribsPower1}.) Bercovici \cite{Bercovici1} subsequently
obtained a new approach to hyper-reflexivity
by means of a very general  result for
spaces with property $X_{\theta,\gamma}$. This property holds for
operator algebras whose commutants contain a pair of commuting
isometries with orthogonal ranges and so is immediately applicable to
free semigroupoid algebras as well as free semigroup algebras.
In the present paper we deal with the case of
finite graphs whose algebras do not have this property and
combining the approaches we obtain
the following main result.
\medskip

\noindent {\bf Theorem} {\it The  free semigroupoid algebra
$\Ll_G$ of a finite directed graph $G$ is hyper-reflexive.}
\medskip

In fact our proof also shows that $\Ll_G \otimes \Bb (\Kk)$ is
hyper-reflexive and we remark on such complete hyper-reflexivity
in Section 3.

Free semigroupoid algebras may be  defined as follows. Let $G$ be
a directed graph with vertex set $\Vv(G)$ and edge set $\Ee(G)$,
each set being finite or countable. A directed path $p$ in $G$ is
either a single vertex (considered as a degenerate path) or an
$N$-tuple of edges $p= e_Ne_{N-1}\dots e_1$ such that the source
vertex of $e_{i+1}$ is equal to the range vertex of $e_i$ for each
$i$.
The length  $|p|$  of the path $p$ is the number of  edges in $p$.
The discrete semigroupoid $\FF_+(G)$ of $G$ is the set of all
finite directed paths with the partially defined associative
product arising from concatenation. Note that $\FF_+(G)$ is a
multiplicatively closed  subset of the (free) path groupoid
$\FF(G)$ of $G$ which is obtained in a similar way from the set of
vertices and all finite paths (or words) in the edges $e$ and
their  formal inverses $e^{-1}$. These words are free except for
the groupoid relations $e^{-1}e = x, ee^{-1} = y$, where $e$ has
$x, y$ as source and range vertices. Also $\FF_+(G)$ is a unital
subset in that it contains the units (vertices) of the  groupoid
$\FF(G)$. These observations and the parallel with free semigroups
support the use here of the term free semigroupoid.

The path semigroupoid $\FF_+(G)$ gives rise in a natural way to a
family  of partial isometries: Let $\Hh$ be the Hilbert space with
basis  $\{\xi_w  :  w\in \ \FF_+(G)\}$ and define the partial
isometry $L_w$ for which
$$ L_w\xi_{w'} =  \left\{%
\begin{array}{ll}
    \xi_{ww'} & \hbox{if } ww' \hbox{ is defined}\\
    0 & \hbox{otherwise.} \\
\end{array}%
\right.
$$
If $x\in \Vv (G)$, then $L_x$ is the projection onto
the closed space spanned by
 $\{\xi_{w} : w = xu\}$. Thus, $L_x$
is the initial projection $L_e^*L_e$ for any $e$ with
$e = ex$.
\begin{defn}
Let $G$ be a countable directed graph. Then the free
semigroupoid algebra $\Ll_G$ is the {\small WOT}-closed
algebra generated by $\{L_e \; :\; e\in \Ee(G)\cup \Vv(G)\}$.
\end{defn}

In particular, a graph with a single vertex and $n$ edges gives
rise to the free semigroup algebra $\Ll_n$.

We would like to thank Ken Davidson for some helpful comments on a
preliminary version of this paper.

\section{\sf Some general results about hyper-reflexivity.}
In this section, we  recall some key results about
hyper-reflexivity. Moreover we derive some auxiliary results that
we need for our proof, including the fact that certain partial
inflations of a hyper-reflexive algebra are again hyper-reflexive.

\begin{thm} (Davidson \cite{Davidson1}.) The algebra $H^\infty(\TT)$ realized as
the algebra of Toeplitz operators acting on $H^2(\TT)$ is
hyper-reflexive, with constant at most $19$.
\end{thm}
Recall that the {\small WOT}-continuous functionals on $\Bb(\Hh)$
are exactly the functionals of the type $\phi= \displaystyle
\sum_{i= 1}^n{[x_i\otimes y_i]}$ where the elementary functionals
$[x \otimes y]$ are defined by their action on an operator $T$ by
$[x \otimes y](T) =  (Tx,y)$.
 A {\small WOT}-closed space
$\Mm$ is said to have property $(\AA_1)$ if every {\small
WOT}-continuous functional $\phi$ on $\Mm$ can be written $\phi=
[x\otimes y]$ for some $x,y\in\Hh$. Moreover, given $r>0$, $\Mm$
is said to have property $(\AA_1(r))$ if for any $\epsilon >0$,
one can find vectors $x,y\in\Hh$ such that $\phi= [x\otimes y]$
and $\|x\|,\; \|y\| \leq
(r+\epsilon)\|\phi\|^{1/2}$.
\begin{thm}  (Davidson \cite{Davidson1}.) Suppose $\Aa$ is a hyper-reflexive
subspace of $\Bb(\Hh)$ with property $(\AA_1 (r))$ for some $r>0$.
Then every {\small WOT}-closed subspace of $\Aa$ is
hyper-reflexive.
\end{thm}

\begin{cor} Every {\small WOT}-closed subspace of $H^\infty(\TT)$,
as an operator algebra on $H^2(\TT)$, is hyper-reflexive.
\end{cor}

We remark that the proof of Theorem 2.1 is rather subtle whereas
the hyper-reflexivity of $H^\infty(\TT)$ as an operator algebra on
$L^2(\TT)$ is a consequence of the hyper-reflexivity of $L^\infty(\TT)$
on $  L^2(\TT)$
and the elementary Theorem 2.2.

The following theorem is due to Bercovici and applies immediately
to many free semigroupoid algebras.

\begin{thm} (Bercovici
\cite{Bercovici1}.) \label{2.4}Suppose that $\Aa$ is a {\small
WOT}-closed subspace of $\Bb(\Hh)$ such that there exist two
isometries in the commutant of $\Aa$ with orthogonal ranges. Then
$\Aa$ is hyper-reflexive, with constant at most $3$.
\end{thm}

The next theorems are immediate consequences of this double
isometry theorem and the fact that the commutant of a free
semigroup algebra, or free semigroupoid algebra, is the
companion algebra determined by the partial isometries  $R_e$ for
which $R_e\xi_u = \xi_{ue}$.

\begin{thm} (Davidson and Pitts \cite{DP1}.)  The free
semigroup algebras $\Ll_n, $ for $ 2 \le n \le \infty$, are
hyper-reflexive.
\end{thm}

\begin{thm} (Kribs and Power \cite{KribsPower1}.) Let $G$ be a directed graph
such that from every vertex $v$ there are at least two distinct
directed cycles based at $v$. Then the free semigroupoid algebra
$\Ll_G$ is hyper-reflexive.
\end{thm}

\noindent {\bf Remark.}
Bercovici obtained
hyper-reflexivity  for a larger class of spaces
than those given in Theorem 2.4, namely those with
 property
$X_{\theta,\gamma}$. This property is defined as follows.
For $x,y \in\Hh$ let
$[x\otimes y]_\Mm$ denote the restriction of the elementary functional to the
space $\Mm$. Given a number $\theta>0$  the set
$\Xx_\theta(\Mm)$ consists of all norm continuous functionals
$\phi$ on $\Mm$ with the following property: for every finite set
$\Ff\subset \Hh$ and every positive number $\epsilon$, there exist
vectors $x,y\in \Hh$ such that:
\begin{enumerate}
    \item $\|x\|,\|y\| \leq 1$;
    \item $\|\phi-[x\otimes y]_\Mm\|\leq \theta+\epsilon$;
    \item $\|[x\otimes f]_\Mm\|+\|[f\otimes y]_\Mm\|<\epsilon$ for
    $f\in \Ff$.
\end{enumerate}
It is well known that $\Xx_\theta (\Mm)$ is norm closed, convex
and balanced. For $\gamma>0$, the space $\Mm$ is said to have
property $X_{\theta,\gamma}$ if the set $\Xx_\theta (\Mm)$
contains all the {\small WOT}-continuous functionals on $\Mm$ with
norm no greater than
$\gamma$.

We now consider spaces of operators with various forms of block matrix structure.

Let us say that a {\small WOT}-closed space of operators $\Ss$ is
a  block matrix space if there exist two families of projections
$\Pp =
\{P_i\; ; i= 1,\dots , r \}$ and $\Qq =  \{Q_i\; ; i= 1,\dots , s \}$
which partition the identity operator, and  a subset
 $\Ff\subset \Qq \times \Pp$ such that
 \[
\Ss= \{T\in\Bb(\Hh)\; : \; QTP= 0,\ \hspace*{5mm}\mbox{for all  }
    (Q,P) \in \Ff\}.
\]
The following proposition is elementary.
\begin{prop} Let $\Ss$ be a block matrix space
defined by the family $\Ff$.
Then $\Ss$ is hyper-reflexive with constant not exceeding
$n$ where $n= |\Ff |$.
\end{prop}

\begin{prop} \label{2.6}Suppose that $\Ss$ is a hyper-reflexive space
and let $S_0$ be a {\small WOT}-closed subspace of $\Ss$.
Then $\Ss_0$ is hyper-reflexive if there exists a constant $C>0$
such that:
$$d(X,\Ss_0)\leq C \beta_{\Ss_0}(X) \hspace*{5mm} \mbox{ for all } X\in \Ss .$$
\end{prop}
\begin{demo}
     For a given space $\Ss$
      write $\Ff_\Ss= \Pp\times \Qq$ for the set of pairs of
      projections annihilating $\Ss$.
Since $\Ss_0 \subset \Ss$, we have $\Ff_\Ss \subset
\Ff_{\Ss_0}$ and so, for $T$ in $\Bb(\Hh)$,
$\beta_\Ss(T)\leq\beta_{\Ss_0}(T)$.
It follows that there exists
$S\in\Ss$ such that $\|T-S\|\leq C_1 \beta_{\Ss_0}(T)$
where $C_1$ is the hyper-reflexivity constant of $\Ss$.
Thus
$$d(T,\Ss_0)\leq \|T-S\| + d(S,\Ss_0)\leq
C_1\beta_{\Ss_0}(T) + d(S,\Ss_0)$$
and so $\Ss _0$ is
hyper-reflexive with constant no greater than $C_1+C$.
\end{demo}

We now give an application of these elementary propositions to
subspaces of block matrix spaces.

\begin{prop} \label{2.7}Let $\Mm$ be a {\small WOT}-closed subspace of
$\Bb(\Hh)$ which is contained in a block matrix space associated
with the family $\Ff\subset \Qq \times \Pp$. Moreover suppose that
that
\begin{enumerate}
    \item $\displaystyle \Mm= \Sigma_{(Q,P)\in\Ff}\oplus Q\Mm P$.
    \item For each $ (Q,P)\in \Gg$ the subspace $ Q\Mm P $ is
hyper-reflexive.
\end{enumerate}
Then the space $\Mm$ is hyper-reflexive.
\end{prop}
\begin{demo} The space $\Mm$ is a subspace of the block
subspace $\displaystyle \Ss = \Sigma_{(Q,P)\in\Ff}\oplus Q\Bb(\Hh)
P$ and so, by the preceding two propositions, it is enough to
prove that $d(T,\Mm)\leq C \beta_\Mm(T)$ for $T\in \Ss$. Since the
entries in the matrix decomposition of $\Mm$ are independent and
$\beta_\Mm(QTP)\leq \beta_\Mm(T)$, it is enough to prove the
inequality for $T\in Q\Bb(\Hh)P$ for any pair $(Q,P)$. But the
hyper-reflexivity of $Q\Mm P$ gives a constant $C>0$ such that
$d(T,Q\Mm P)\leq C \beta_{Q\Mm P}(T)$. So $d(T,\Mm) = d(T,Q\Mm
P)\leq C \beta_{Q\Mm P}(T)\leq C\beta_\Mm(T)$, where the last
inequality holds since the annihilating set $\Ff_{Q\Mm P}$ for
$Q\Mm P$ contains $ \Ff_\Mm$.
\end{demo}

The following corollary now follows from Davidson's theorem for
$H^\infty(\TT)$.
\begin{cor} Let $\Mm\subseteq M_n(\CC)$ be a bimodule for the space
$D_n$ of diagonal matrices. Then the space $H^\infty(\TT)\otimes
\Mm$, as a space of operators on $H^\infty(\TT)\otimes \CC^n$, is
hyper-reflexive.
\end{cor}
We now turn to matrix spaces and algebras where the matrix entries
are not independent. We start with a straightforward application
of Bercovici's theorem to give a short new proof of a known fact
(see \cite{dav-book}).

\begin{prop} Let $\Aa$ be a {\small WOT}-closed subspace of $\Bb(\Hh)$.
Then the infinite ampliation $\Aa^{(\infty)}$ is a
hyper-reflexive subspace of $\Bb(\Hh^{(\infty)})$.
\end{prop}
\begin{demo}
The space $\Aa^{(\infty)}$ is equal to $ \Aa \otimes \CC I$ acting
on $\Hh \otimes \ell^2(\ZZ)$ and so the commutant contains all
operators of the form $I \otimes W$ with $W$ in $\Bb(\ell^2(\ZZ)$.
In particular there are two isometries in the commutant  with
orthogonal ranges and so Bercovici's theorem applies.
\end{demo}\par

The next proposition concerns  what might be termed a partial
inflations of operator spaces.
\begin{prop} \label{2.10} Let $\Aa$ be a  hyper-reflexive
subspace of $\Bb(\Hh)$ and let $\Mm$ be an invariant subspace for
$\Aa$. Then, for any given positive integer $n$, the space
$$
\tilde{\Aa} =  \{A\oplus A_{|\Mm}^{(n)}~;~~ A\in\Aa\}
$$
 is hyper-reflexive.
\end{prop}
\begin{demo} We start with $n= 1$.
Since the space $\Ss =\Bb(\Hh) \oplus \Bb(\Mm)$ is a
hyper-reflexive subspace of $\Bb(\Hh \oplus\Mm)$ it follows from
Proposition
\ref{2.6} that  it is sufficient  to show that there exists a
constant $C$ such that $d(\tilde{X},\tilde{\Aa})\leq C
\beta_{\tilde{\Aa}}(\tilde{X}) \mbox{ where } \tilde{X}= X\oplus
Y, \; X\in \Bb(\Hh) \mbox{ and } Y\in \Bb(\Mm).$ We consider the
Hilbert space decomposition $\tilde{\Hh}= \Hh\oplus \Mm =
(\Hh_1\oplus\Mm)\oplus\Mm$ and the projection defined by the
matrix
$$P= \frac{1}{2}
\begin{pmatrix}
  0 & 0 & 0 \\
  0 & I & I \\
  0 & I & I \\
\end{pmatrix}.
$$
Relative to the decomposition $\Hh_1 \oplus \Mm$   write
$$X=
\begin{pmatrix}
  X_1 & 0 \\
  X_2 & X_3 \\
\end{pmatrix},
$$ so that
$$P^\perp\tilde{X}P =  \frac{1}{4}\begin{pmatrix}
  0 & 0 & 0 \\
  0 & X_3-Y & X_3-Y \\
  0 & Y-X_3 & Y-X_3 \\
\end{pmatrix} \mbox{~~~~~~~~and~~~~~~~}
\beta_{\tilde{\Aa}}(\tilde{X})\geq
\frac{1}{4}
\|X_3-Y\|.$$
 Now write
$$\tilde{X}= \begin{pmatrix}
  X & 0 \\
  0 & X_3 \\
\end{pmatrix} + \begin{pmatrix}
  0 & 0 \\
  0 & Y-X_3 \\
\end{pmatrix} \mbox{ ~~~~~~~and set~~~~~~~ } X'= \begin{pmatrix}
  X & 0 \\
  0 & X_3 \\
\end{pmatrix}.
$$
Then $$d(\tilde{X},\tilde{\Aa})\leq
d(X',\tilde{\Aa})+4\beta_{\tilde{\Aa}}(\tilde{X}).$$
Also $\|X'-A\oplus A_{|\Mm}\|= \|X-A\|$. Since $\Aa$ is
hyper-reflexive with constant $C_\Aa$, we have
$$d(X',\tilde{\Aa})= d(X,\Aa)\leq C_\Aa\beta_\Aa(X) \leq C_\Aa
\beta_{\tilde{\Aa}}(\tilde{X}).$$ The last inequality holds since
$\Ff_{\tilde{\Aa}} \supset \Ff_\Aa$ and $X$ is a summand of $\tilde{X}$.
Finally we have
$$d(\tilde{X},\tilde{\Aa})\leq
(C_\Aa+4)\beta_{\tilde{\Aa}}(\tilde{X})$$
 as desired.

Suppose now  that $\tilde{\Aa}_n = \{A\oplus
A_{|\Mm}^{(n)} \; ; \;A\in\Aa\}$ is hyper-reflexive. Since
$\tilde{A}_{n+1}= A\oplus A_{|\Mm}^{(n+1)}$ is unitarily  equivalent
to $\tilde{A}_n\oplus (\tilde{A}_{n})_{|0\oplus\Mm}$, it follows
that $\tilde{\Aa}_{n+1}$ is also
hyper-reflexive.
\end{demo}

The next  proposition deals with infinite partial
ampliations:

\begin{prop} \label{2.11}Let $\Aa$ be a {\small
WOT}-closed operator algebra which is hyper-reflexive and  let $\Mm$ be an invariant
subspace for $\Aa$. Then the
operator algebra
$$\tilde{\Aa} =  \{A\oplus
A_{|\Mm}^{(\infty)} \; ;\; A\in\Aa\}$$
is hyper-reflexive.
\end{prop}
\begin{demo}
Since $\Bb(\Hh)\oplus\Aa_{|\Mm}^{(\infty)}$ is hyper-reflexive, as
a simple consequence of Proposition 2.11, it is enough in view of
Proposition 2.8 to derive the distance formula for elements of the
form $X\oplus A_{|\Mm}^{(\infty)}$, $X\in \Bb(\Hh)$, $A\in \Aa$.
Set $\tilde{\Aa}_1=  \{A\oplus A_{|\Mm} ; A\in \Aa\}$, which is
hyper-reflexive by the previous proposition, with constant $C_1$.
Then
$$d(X\oplus A_{|\Mm}^{(\infty)}, \tilde{\Aa})= d(X\oplus
A_{|\Mm}, \tilde{\Aa}_1)\leq C_1\beta_{\tilde{\Aa}_1}(X\oplus
A_{|\Mm})\leq C_1\beta_{\tilde{\Aa}}(X\oplus
A_{|\Mm}^{(\infty)})$$ which completes the proof.
\end{demo}


\section{\sf Hyper-reflexivity of free semi{groupo\"{\i}}d algebras}
We now consider a finite directed graph $G$ with semigroupoid
algebra $\Ll_G$ acting on the Fock space $\Hh_G$ and, as we have
intimated in the introduction, we shall show that $\Ll_G \otimes
\Bb(\Kk)$ is hyper-reflexive, for $\Kk$ a separable Hilbert space.

It does not seem to be known whether the hyper-reflexivity of an
operator algebra $\Aa$ entails that of $\Aa \otimes \Bb(\Kk)$,
even allowing for a different hyper-reflexivity constant. We find
it convenient in the subsequent discussion to say that $\Aa$ is
{\it completely hyper-reflexive} if this apparently stronger
assertion holds. That completely hyper-reflexive spaces are
hyper-reflexive is elementary.

Our proof uses an induction argument with an induction step that
describes a block matrix structure of $\Ll_G$ in terms of
$\Ll_{G'}$ where $G'$ is a subgraph of $G$ obtained by deleting a
single edge $e$. We can in fact restrict attention to the case
where $e$ does not lie on any directed loop path with source
vertex equal to the source of $e$. In the block matrix
decomposition there is column subspace of the form $\Ss \otimes
\ell^2$ with $\Ss \subseteq \Ll_{G'}$ a hyper-reflexive space, and
so we need complete hyper-reflexivity for $\Ll_{G'}$ even to
deduce that $\Ll_G$ is hyper-reflexive.

We shall make use of the following strengthened form of Davidson's
theorem for $H^\infty(\TT)$ on $H^2(\TT)$ the proof of which is an
entirely routine adaptation of the proof given in
\cite{Davidson1}.

\begin{thm}  The algebra $H^\infty(\TT)$ realized as
the algebra of T{\oe}plitz operators acting on $H^2(\TT)$ is
completely hyper-reflexive, with constant at most $19$.
\end{thm}

We first  isolate the case when $G$ is a cycle; in this case we
write $G= C_n$, where $n$ is the number of vertices in $G$. The
corresponding Fock space will be denoted by $\Hh_n$.
\begin{thm} \label{3.1} For any integer $n>0$, the cycle algebra
$\Ll_{C_n}$ is completely hyper-reflexive.
\end{thm}
\begin{demo} If $n= 1$,
then  $\Ll_{C_1} = H^\infty(\TT)$, acting on $H^2$, which is
completely hyper-reflexive. For $n\geq 2$, the Fock space $\Hh_n$
decomposes as an orthogonal sum of the ranges of the projections
$L_v$ for $v \in \Vv(C_n)$:
\[
\displaystyle \Hh_n= \bigoplus_{i= 1}^n L_{v_i}\Hh.
\]
Identify each space $L_{v_i}\Hh$  with the Hardy space $H^2(\TT)$.
Then it is straightforward to verify that
 $\Ll_{C_n}$ is unitarily
equivalent to the matrix function algebra
acting on $H^2\oplus \cdots \oplus H^2$, with elements of the form
$$\begin{pmatrix}
  f_{11}(z^n) & z^{n-1}f_{12}(z^n) & z^{n-2}f_{13}(z^n) & \cdots & z
  f_{1n}(z^n) \\
  zf_{21}(z^n) & f_{22}(z^n) & z^{n-1}f_{23}(z^n) & \cdots & z^2f_{2n}(z^n) \\
  z^2f_{31}(z^n) & zf_{32}(z^n) & f_{33}(z^n) & \cdots & \vdots \\
  \vdots & \vdots & \vdots & \ddots & \vdots \\
  z^{n-1}f_{n1}(z^n) & z^{n-2}f_{n2}(z^n) & \cdots & \cdots & f_{nn}(z^n) \\
\end{pmatrix},$$
where each $f_{ij} \in H^\infty(\TT)$.
 Thus the algebra $\Ll_{C_n}$ is
identifiable with a block operator matrix algebra with copies of
$H^\infty(z^n)$ on the diagonal and with off-diagonal spaces of
the form $z^jH^\infty(z^n)$, where $H^\infty(z^n)$ denotes the
subalgebra spanned by the powers $z^{nk}, k=0,1,2,\dots .$


The space $H^\infty(z^n)$ is unitarily equivalent to the
ampliation algebra $H^\infty(\TT)^{(n)} = H^\infty(\TT)\otimes \CC
I$ acting on $H^2(\TT) \otimes \CC^n$. Also it is elementary to
see that finite  ampliations of a hyper-reflexive algebra are
hyper-reflexive. Thus $H^\infty(z^n)$ on $H^2(\TT)$ is completely
hyper-reflexive. Also it is elementary that if $\Ss$ is a
completely hyper-reflexive space and $U$ is an isometry (such as
$T_{z^j} \otimes I$ in the present context) then $U\Ss$ is a
completely hyper-reflexive space. Thus $z^jH^\infty(z^n)$ on
$H^2(\TT)$ is completely hyper-reflexive for each $j = 1,\dots,n$.
It follows now, from Proposition 2.9 for examples, that the matrix
function operator algebra above is completely hyper-reflexive, as
required.
\end{demo}

\medskip

We now turn to the following key induction step.

\begin{prop} \label{3.2}Let $G$ be a graph with directed edge $e= yex$,
for which there is no directed path from vertex $y$ to vertex $x$
and let $G'$ be the graph where $e$ has been deleted : $G'=
G\backslash e$. If $\Ll_{G'}$ is completely hyper-reflexive, then
so too is $\Ll_G$.
\end{prop}
\begin{demo} We write $\Hh= \Hh_G$ and $\Hh'= \Hh_{G'}$ for
the Fock spaces associated with $G$ and $G'$ respectively. Let
$\Hh'$ be viewed as the natural subspace of $\Hh$, that is, the
subspace spanned by the basis elements $\xi _w$ for which $e$ does
not appear in $w$. We  write $e\not \in w$ to indicate this.
With respect to the decomposition $\Hh= \Hh' \oplus \Kk$ note
that $\Kk \subseteq \mbox{ ker}L_e$ and that the operator matrix
for $L_e$ is strictly lower triangular. Also, for any other edge
$g\in\Ee (G)$, $L_g$ has block diagonal matrix. It follows that
for a given word $w\in \FF_+(G)$, $L_w$ is block diagonal if
$e\not \in w$, and is lower triangular if $e\in w$. Thus $\Ll_G$
is the direct sum of its diagonal subalgebra, which we denote by
$\Dd$, and its lower triangular subalgebra, $\Tt$ say.

We  now show that both $\Dd$ and $\Tt$ are completely
hyper-reflexive and it follows readily (by a simple adaptation of
the proof of Proposition 2.9 for example) that $\Dd + \Tt$ is
completely hyper-reflexive.

Consider first the decomposition of $\Kk$ associated with the
different words in $G'$ ending at $x$:
$$\Kk= \bigoplus _{i= 1}^n \Kk^i, \hspace*{5mm} \Kk^i\simeq\Hh_y, $$
where $\Hh_y= \mbox{span}\{\xi_w; \hspace*{3mm} w= wy\}$ and $n$
is the cardinality of
$\{w\in \Ww (G'); \hspace*{3mm} w= xw\}$. (Possibly $n=
\infty$.) The block diagonal algebra $\Dd$ is therefore unitarily
equivalent to the algebra
$$\{A\oplus \bigoplus_{i= 1}^n A_{|\Kk^i}\hspace*{5mm} A\in
\Ll_{G'}\}.$$ In view of the complete hyper-reflexivity of
$\Ll_{G'}$ it follows from  Propositions \ref{2.10} and \ref{2.11}
that this algebra is completely hyper-reflexive.

To see the complete hyper-reflexivity of the lower triangular
subalgebra $\Tt$, we use a different  orthogonal decomposition of
$\Kk$ namely the one obtained from  the wandering subspace
$$\Kk_0= L_e\Hh \;(= L_eL_e^*\Hh \simeq
L_x\Hh = L_x\Hh'),$$ that is:
$$\Kk= \bigoplus_{w= wy}L_w\Kk_0.$$
With respect to the associated block matrix structure and the
natural identification of each summand  $L_w\Kk_0$ with $\Kk_0$
the restriction operator ${L_{w}}_{|\Kk}$ is a block
matrix whose entries are either the identity or zero.
Now, with respect to the decomposition
\[
\Hh= (\Hh' \ominus
L_x\Hh')\oplus L_x\Hh'\oplus \Kk_0 \oplus(\Kk \ominus \Kk_0),
\]
the operator $L_e$ has the  operator matrix
\[
L_e=
\begin{pmatrix}
  0 & 0 & 0 & 0 \\
  0 & 0 & 0 & 0 \\
  0 & I & 0 & 0 \\
  0 & 0& 0 & 0 \\
\end{pmatrix}.
\]
Let us write the restricted algebra $\Dd_{|\Hh'}$ which is
isomorphic to $\Ll_{G'}$ in terms of its  $2\times2$ block structure  as
$$\Dd_{|\Hh'}= \begin{pmatrix}
  \Aa_1 & \Aa_2 \\
  \Bb & \Cc \\
\end{pmatrix}.$$
Since $\Ll_{G'}$ is completely hyper-reflexive, by assumption, it
follows that the space of matrices
$$L_x\Dd =\begin{pmatrix}
  0 & 0 \\
  \Bb & \Cc \\
\end{pmatrix}$$
is completely hyper-reflexive.  Also $\Tt$ can be identified with
the lower triangular space
$$\Tt \simeq \begin{pmatrix}
  0 & 0 \\
  [\Bb ~~~ \Cc]\otimes l_+^2(E) & 0 \\
\end{pmatrix},$$
where $E= \Ww (G)y$ is the set of words $w$ in $G$ of the form $w=
wy$, and so it follows that  $\Tt$ is completely hyper-reflexive,
as required.
\end{demo}

\begin{thm} For a finite directed graph $G$ the free semigroupoid
algebra $\Ll_G$ is completely hyper-reflexive.
\end{thm}

\begin{demo} If $G$ is the single cycle graph $C_n$ then $\Ll_G$ is completely
hyper-reflexive as we saw above, so we may assume that $G$ is not
of this form.

Suppose first that $G$ is a transitive directed graph, that is,
that every pair of vertices lies on a directed cycle. The
commutant of $\Ll_G \otimes \Bb (\Kk)$ contains the ampliated
commutant $(\Ll_G)'\otimes \mathbb{C} I$. Also, $(\Ll_G)'$ is
isomorphic to $\Ll_{G^t}$ where the transpose graph $G^t$ is also
transitive. (See \cite{KribsPower1}.) In view of our initial
assumption every vertex of the transitive graph, and its
transpose, lies on two distinct cycles. An elementary argument,
given in \cite{KribsPower1}, shows that in this case there is a
pair of isometries in $\Ll_{G^t}$  with orthogonal ranges. Thus
there is a pair of isometries in the commutant of $\Ll_{G}\otimes
\Bb (\Kk)$ with orthogonal ranges and Bercovici's theorem shows
that $\Ll_G$ is completely hyper-reflexive.

Suppose now, by way of induction, that $\Ll_G$ is completely
hyper-reflexive for all graphs with $n $ or fewer edges, and let
$G$ have $n+1$ edges. If $G$ is transitive then it is completely
hyper-reflexive, by the above, and so we may assume that there is
an edge $e = (x,y)$ for which there is no path from $x$ to $y$.
The induction hypothesis and Proposition 3.3 complete the proof.
\end{demo}

Fr\'{e}d\'{e}ric Ja\"{e}ck : {jaeck@math.u-bordeaux1.fr}

Stephen~C. Power : {s.power@lancaster.ac.uk}

\end{document}